\documentclass[a4paper,11pt]{amsart}
\usepackage{amssymb,mathtools}
\usepackage{microtype}
\usepackage{xcolor}
\definecolor{citeblue}{RGB}{0,0,128}
\definecolor{refred}{RGB}{128,0,0}
\usepackage[colorlinks=true,citecolor=citeblue,linkcolor=refred,urlcolor=citeblue]{hyperref}
\hypersetup{pdftitle={Standard graded nearly Gorenstein Ehrhart rings have unimodal h-vectors},pdfauthor={Takayuki Hibi and Sora Miyashita},pdfsubject={Canonical trace ideals, Frobenius splitting, Ehrhart rings, and h-vectors},pdfkeywords={Nearly Gorenstein ring, Ehrhart ring, canonical trace ideal, Frobenius splitting, strong Lefschetz property, graded domain, unimodality, level ring}}

\DeclareMathOperator{\tr}{tr}

\DeclareMathOperator{\Hilb}{Hilb}

\DeclareMathOperator{\indeg}{indeg}
\DeclareMathOperator{\relint}{relint}
\DeclareMathOperator{\Hom}{Hom}
\DeclareMathOperator{\height}{ht}

\newcommand{\kk}{k}
\newcommand{\QQ}{\mathbb{Q}}
\newcommand{\ZZ}{\mathbb{Z}}
\newcommand{\RR}{\mathbb{R}}
\newcommand{\m}{\mathfrak{m}}
\newcommand{\natur}{\natural}

\newtheorem{theorem}{Theorem}[section]
\newtheorem{lemma}[theorem]{Lemma}
\newtheorem{proposition}[theorem]{Proposition}
\newtheorem{corollary}[theorem]{Corollary}
\newtheorem{conjecture}[theorem]{Conjecture}

\newtheorem{maintheorem}{Theorem}

\theoremstyle{definition}
\newtheorem{example}[theorem]{Example}
\newtheorem{remark}[theorem]{Remark}

\title[Unimodality for nearly Gorenstein Ehrhart rings]{Standard graded nearly Gorenstein Ehrhart rings have unimodal $h$-vectors}
\author{Takayuki Hibi}
\author{Sora Miyashita}
\address{Department of Pure and Applied Mathematics, Graduate School of Information Science and Technology, The University of Osaka, Suita, Osaka 565-0871, Japan}
\email[T.H.]{hibi@math.sci.osaka-u.ac.jp}
\email[S.M.]{u804642k@ecs.osaka-u.ac.jp}
\date{November 21, 2026}
\subjclass[2020]{Primary 13H10, 05E40; Secondary 13A02, 52B20}
\keywords{Nearly Gorenstein ring, Ehrhart ring, canonical trace ideal, Frobenius splitting, strong Lefschetz property, graded domain, unimodality, level ring}

\begin{document}

\begin{abstract}
We prove that the $h$-vector $(h_0,\ldots,h_s)$ of every standard graded nearly Gorenstein Ehrhart ring is unimodal, with a peak at position $\lceil s/2\rceil$. More generally, a standard graded Ehrhart ring has an alternatingly increasing $h$-vector if the radical of the ideal generated by the degree-one part of its canonical trace ideal contains the homogeneous maximal ideal. In characteristic zero or two, the generic Artinian reduction of each such ring has the strong Lefschetz property. In characteristic two, both conclusions extend to all $F$-pure standard graded Cohen--Macaulay algebras satisfying the same trace condition. We also prove unimodality for every standard graded nearly Gorenstein domain of dimension at least two with $h$-vector $(1,h_1,\ldots,h_s)$, where $s\geq2$, provided that $h_1-h_s\leq1$. If the ring is level, its $h$-vector is both concave and alternatingly increasing. When $h_1=h_s$, the ring is necessarily level and its $h$-vector has the form $(1,h_1,\ldots,h_1)$.
\end{abstract}

\maketitle

\section{Introduction}

Let $R=\bigoplus_{j\geq0} R_j$ be a standard graded Cohen--Macaulay algebra over a field $\kk$. Set $d=\dim R$ and let $\m_R=\bigoplus_{j>0} R_j$ be its \emph{homogeneous maximal ideal}. The \emph{Hilbert series} $\Hilb_R(t)=\sum_{j\geq0}(\dim_{\kk} R_j)t^j$ has the form $h_R(t)/(1-t)^d$, where $h_R(t)=1+h_1t+\cdots+h_st^s$ and $h_s\neq0$. We call $h_R(t)$ the \emph{$h$-polynomial} and $h(R)=(1,h_1,\ldots,h_s)$ the \emph{$h$-vector} of $R$. The $h$-vector is \emph{unimodal} if $h_0\leq\cdots\leq h_j\geq\cdots\geq h_s$ for some $j$ with $0\leq j\leq s$. It is \emph{log-concave} if $h_i^2\geq h_{i-1}h_{i+1}$ for every $i$ with $1\leq i<s$.

By classical results of Stanley, the $h$-vector of every standard graded Gorenstein algebra is symmetric \cite[Theorem~4.1]{S78}, and every standard graded Cohen--Macaulay domain with a symmetric $h$-vector is Gorenstein \cite[Theorem~4.4]{S78}. A symmetric $h$-vector need not be unimodal. For example, Stanley constructed an Artinian Gorenstein algebra with $h$-vector $(1,13,12,13,1)$ \cite[Example~4.3]{S78}.

The following conjecture goes back to Hibi and Stanley and remains open.

\begin{conjecture}[{\cite[Conjecture~1.5]{H89}, \cite[Conjecture~4~(a), p.~512]{S89}}]\label{conj:classical}
Every standard graded Gorenstein domain has a unimodal $h$-vector.
\end{conjecture}

Stanley's original conjecture asked for the stronger property of log-concavity for standard graded Cohen--Macaulay domains, with the Gorenstein case suggested as a possible restriction. Conca, Singh, and Varbaro constructed semi-standard graded Gorenstein normal domains with non-unimodal $h$-vectors \cite[Theorem~4.2]{CSV25}. Their examples show that the standard grading hypothesis in Conjecture~\ref{conj:classical} is essential.

Let $P\subset\RR^q$ be a lattice polytope. Its \emph{Ehrhart ring} $A_P=\kk[\mathbf{x}^{\boldsymbol u}t^n:\boldsymbol u\in nP\cap\ZZ^q,\ n\geq0]$, graded by the exponent of $t$, is a normal Cohen--Macaulay domain \cite[Theorem~1]{Hoc72}. It is standard graded precisely when $P$ has the \emph{integer decomposition property} (IDP), meaning that for every integer $n\geq1$, each lattice point of $nP$ is a sum of $n$ lattice points of $P$. The Hilbert series of $A_P$ has the form $h_P^*(t)/(1-t)^{\dim P+1}$. The numerator is the \emph{$h^*$-polynomial} $h_P^*(t)=\sum_{i=0}^s h_it^i$, where $h_s\ne0$. Its coefficient vector $h^*(P)=(h_0,\ldots,h_s)$ is the \emph{$h^*$-vector} of $P$. If $P$ is IDP, then $h^*(P)=h(A_P)$. Let $\relint$ denote \emph{relative interior}, taken in the affine span. We call $P$ \emph{Gorenstein} if $A_P$ is Gorenstein.

In 1992 Hibi conjectured unimodality of $h^*(P)$ for every reflexive lattice polytope \cite[p.~111]{H92}. Here a full-dimensional lattice polytope with the origin in its interior is \emph{reflexive} if its polar is a lattice polytope. Musta\c{t}\u{a} and Payne gave a counterexample \cite[Example~1.1]{MP05}, and Payne constructed further counterexamples among reflexive simplices \cite[Theorem~1.4]{P08}. These examples fail IDP. For standard graded Gorenstein toric rings, unimodality was also conjectured by Ohsugi and Hibi \cite[Introduction]{OH06}.

Bruns and R\"omer proved unimodality for Gorenstein polytopes with a regular unimodular triangulation \cite[Theorem~1]{BR07}. Adiprasito, Papadakis, Petrotou, and Steinmeyer announced results removing the triangulation hypothesis for Gorenstein IDP polytopes \cite{APPS22}. Adiprasito, Papadakis, and Petrotou proved unimodality for reflexive IDP polytopes \cite[Theorem~1.4]{APP25} (see also \cite[Corollary~6.3]{AKOPP26}). By \cite[Corollary~7]{BR07}, every Gorenstein IDP polytope has the same $h^*$-polynomial as a reflexive IDP polytope, so unimodality holds for every Gorenstein IDP polytope. On the other hand, Ferroni constructed smooth IDP polytopes with non-unimodal $h^*$-vectors and Gorenstein IDP polytopes with non-log-concave $h^*$-vectors \cite{F26}.

This raises the question whether the $h$-vector remains unimodal when a standard graded Ehrhart ring is nearly Gorenstein. Let $\omega_R$ denote the \emph{graded canonical module} of $R$. The \emph{canonical trace ideal} is $\tr_R(\omega_R)=\sum_{\varphi\in\Hom_R(\omega_R,R)}\varphi(\omega_R)$, where the sum runs over homogeneous maps of all degrees. The ring $R$ is Gorenstein if and only if $\tr_R(\omega_R)=R$ \cite[Lemma~2.1]{HHS19}. A Cohen--Macaulay ring $R$ is \emph{nearly Gorenstein} if $\tr_R(\omega_R)\supseteq\m_R$ \cite{HHS19}. The ring $R$ is \emph{level} if $\omega_R$ is generated in one degree and \emph{anticanonical level} if $\Hom_R(\omega_R,R)$ is generated in one degree (cf.\ \cite[Definition~3.3]{Pa19}).

For IDP $(0,1)$-polytopes $P$, the nearly Gorenstein property of $A_P$ was characterized in \cite[Theorem~34]{HKMM23}.

We call an $h$-vector \emph{alternatingly increasing} if its entries are nondecreasing when read alternately from the two ends, starting on the left.

We use the following condition on the canonical trace ideal \cite[Definition~1.1]{M26}:
\[
(\natur)\qquad \sqrt{[\tr_R(\omega_R)]_1R}\supseteq\m_R.
\]
The following is our first main result.

\begin{maintheorem}\label{thm:A}
Let $P$ be an IDP lattice polytope and let $R=A_P$ have $h$-vector $h(R)=(h_0,\ldots,h_s)$. Suppose that one of the following conditions holds:
\begin{enumerate}
\item $R$ satisfies $(\natur)$;
\item $R$ is level and anticanonical level, and $[\tr_R(\omega_R)]_1\neq0$;
\item $R$ is level and $P$ has an interior lattice point.
\end{enumerate}
Then $h(R)$ is alternatingly increasing, that is,
\[
 h_0\leq h_s\leq h_1\leq h_{s-1}\leq\cdots\leq h_{\lceil s/2\rceil}.
\]
In particular, $h_0\leq\cdots\leq h_{\lceil s/2\rceil}\geq\cdots\geq h_s$, so the $h$-vector is unimodal, with a peak at position $\lceil s/2\rceil$.
\end{maintheorem}

If $R$ is nearly Gorenstein, then $[\tr_R(\omega_R)]_1=R_1$. Thus $(\natur)$ holds, and Theorem~\ref{thm:A}~\textup{(1)} gives the following corollary.

\begin{corollary}\label{cor:ng-ehrhart}
The $h$-vector $(h_0,\ldots,h_s)$ of a standard graded nearly Gorenstein Ehrhart ring is unimodal, with a peak at position $\lceil s/2\rceil$.
\end{corollary}

For an Ehrhart ring, the monomials of the canonical ideal and of its inverse are determined by the lattice points of the cone over $P$. In particular, the hypotheses on the canonical trace ideal in Theorem~\ref{thm:A} are preserved on passing to characteristic two. Since the $h$-vector also depends only on $P$, the resulting inequalities hold over every field.

In characteristic zero or two, the generic Artinian reduction of $R$ in Theorem~\ref{thm:A} has the strong Lefschetz property (Corollary~\ref{cor:ehrhart-consequences}\,\textup{(\ref{item:ehrhart-lefschetz})}). In the Gorenstein case this follows from \cite[Theorem~6.1]{APP25} and the openness statement in \cite[Proposition~13.6]{PP20}.

For level Ehrhart rings, Theorem~\ref{thm:A}~\textup{(3)} answers the unimodality question in \cite{Hi19} and also follows from \cite{APP25} (see Remark~\ref{rem:app-scope}).

In characteristic two, the Ehrhart ring hypothesis in Theorem~\ref{thm:A}~\textup{(1)} can be replaced by $F$-purity. The following is our second main result.

\begin{maintheorem}\label{thm:B}
Let $R$ be an $F$-pure standard graded Cohen--Macaulay algebra over a field $\kk$ of characteristic two. If $R$ satisfies $(\natur)$, then $h(R)$ is alternatingly increasing and the generic Artinian reduction of $R$ has the strong Lefschetz property.
\end{maintheorem}

In particular, Theorem~\ref{thm:B} applies to every $F$-pure standard graded nearly Gorenstein algebra in characteristic two. Its conclusions also hold for every Veronese subring of $R$ (Corollary~\ref{cor:fpure-veronese}).

For a standard graded nearly Gorenstein domain of dimension $d\ge2$, $h_s=1$ forces $R$ to be Gorenstein \cite[Theorem~1.2]{MPG26}. We compare $h_s$ with the \emph{codimension} $h_1=\dim_{\kk}R_1-d$. An $h$-vector is \emph{concave} if $h_i-h_{i-1}\geq h_{i+1}-h_i$ for every $1\leq i<s$.

The following is our third main result.

\begin{maintheorem}\label{thm:C}
Let $R$ be a standard graded Cohen--Macaulay domain over a field, with $\dim R\geq2$ and $h$-vector $h(R)=(1,h_1,\ldots,h_s)$, where $s\geq2$. Suppose that $R$ satisfies $(\natur)$ and that $h_1-h_s\leq1$.
Then $h_s\leq h_1$, and $h(R)$ is unimodal with $\max_i h_i-h_{\lceil s/2\rceil}\leq1$.
In characteristic zero, the generic Artinian reduction of $R$ has the strong Lefschetz property. If $R$ is level, then $h(R)$ is concave and alternatingly increasing. In particular, it is log-concave and has a peak at position $\lceil s/2\rceil$. If $h_1=h_s$, then $R$ is level and $h(R)=(1,h_1,\ldots,h_1)$.
\end{maintheorem}

\vbox{
Theorems~\ref{thm:A}, \ref{thm:B}, and~\ref{thm:C} suggest the following extension of Conjecture~\ref{conj:classical}.

\begin{conjecture}\label{conj:unimodality}
Every standard graded nearly Gorenstein domain has a unimodal $h$-vector.
\end{conjecture}
}

The remainder of the paper is organized as follows.
Section~\ref{sec:canonical} collects the facts about canonical trace ideals and Artinian reductions used below. In Sections~\ref{sec:frobenius} and~\ref{sec:fpure}, we use Frobenius splittings to prove Theorem~\ref{thm:A} for Ehrhart rings and Theorem~\ref{thm:B} for $F$-pure rings, respectively. In Section~\ref{sec:domains}, we prove Theorem~\ref{thm:C} using Veronese subrings of a polynomial ring in two variables.

\section{Canonical trace ideals and Artinian reductions}\label{sec:canonical}

In this section, we describe the degree-one part of the canonical trace ideal and the multiplication pairing on an Artinian reduction. Let $R$ be a $d$-dimensional standard graded Cohen--Macaulay algebra over $\kk$, and set $s=\deg h_R(t)$. We use the grading convention $[M(a)]_j=M_{j+a}$. The \emph{initial degree} of a nonzero finitely generated graded module $M$ is $\indeg M=\min\{j:M_j\ne0\}$.

Let $Q^*(R)$ denote the \emph{graded total ring of fractions} of $R$, obtained by inverting all homogeneous nonzerodivisors. A \emph{graded fractional ideal} is a finitely generated graded $R$-submodule $I\subseteq Q^*(R)$. Set $I^{-1}=(R:I)=\{v\in Q^*(R):vI\subseteq R\}$. If $I$ contains a homogeneous unit of $Q^*(R)$, localization identifies homogeneous maps $I\to R$ with multiplication by homogeneous elements of $I^{-1}$, so $\tr_R(I)=II^{-1}$ \cite[Lemma~2.6]{M26}. A homogeneous ideal is a \emph{canonical ideal} if it is isomorphic, up to a graded shift, to $\omega_R$. Over a domain, a shift of $\omega_R$ can be identified with a graded fractional ideal \cite[Proposition~3.3.18]{BH98}. More generally, if $R$ is reduced and a graded ideal $I\subseteq R$ represents a shift of $\omega_R$, then $I$ contains a homogeneous nonzerodivisor.

In degree one, the product $II^{-1}$ can receive contributions from several degrees of $I$. The next lemma gives conditions under which the lowest-degree component of $I$ suffices.

\begin{lemma}\label{lem:linear-trace-product}
Let $I$ be a graded fractional ideal isomorphic to a shift of $\omega_R$, and suppose that $I$ contains a homogeneous unit of $Q^*(R)$. Set $a=\indeg I$. If $R$ is a domain, $R$ is level, or $\height([\tr_R(\omega_R)]_1R)\ge2$, then
\begin{equation}\label{eq:linear-trace-product}
 [\tr_R(\omega_R)]_1=I_a[I^{-1}]_{1-a}.
\end{equation}
Here the product on the right is the $\kk$-span of the pairwise products.
\end{lemma}
\begin{proof}
We may assume that $\kk$ is infinite. Indeed, both sides of~\eqref{eq:linear-trace-product} commute with extension to $\kk(t)$, since $I$ is finitely presented and $I^{-1}=\Hom_R(I,R)$. Equality descends by faithful flatness.

We have $\tr_R(\omega_R)=II^{-1}$. If $I=uR$, then $I_a=\kk u$ and $I^{-1}=u^{-1}R$, giving~\eqref{eq:linear-trace-product}. This includes $d=0$, when $I=Q^*(R)=R$. We may therefore assume that $d>0$ and $II^{-1}\subseteq\m_R$. If $[II^{-1}]_1=0$, the formula is immediate. Otherwise the trace has initial degree one.

Suppose $R$ is a domain or is level. Since $I$ becomes free of rank one at every minimal prime, $R$ is generically Gorenstein. Hence $I_a$ contains an element with zero annihilator by \cite[Remark~3.14]{M26}, and \cite[Lemma~2.12~(2)]{M26} gives~\eqref{eq:linear-trace-product} in the initial degree of the trace.

It remains to consider $\height([II^{-1}]_1R)\ge2$. Set $X=I_a[I^{-1}]_{1-a}$ and $Y=\sum_{j>a}I_j[I^{-1}]_{1-j}$, so that $[II^{-1}]_1=X+Y$. For $j>a$, the product $I_a[I^{-1}]_{1-j}$ lies in degree $a+1-j\le0$ of the proper ideal $II^{-1}$, so it vanishes. Hence $YI_a=0$.

Since $R$ is Cohen--Macaulay and the field is infinite, we can choose a regular sequence $\theta_1,\theta_2$ in $[II^{-1}]_1$. Since $YI_a=0$, we have $XY=0$. Choose $x_i\in X$ and $y_i\in Y$ with $\theta_i=x_i+y_i$. Then $\theta_2x_1=x_2x_1=\theta_1x_2$. Since $\theta_2$ is regular modulo $\theta_1$, we have $x_1\in\theta_1R$. Both are linear forms, so $x_1=\lambda\theta_1$ for some scalar $\lambda$. If $\lambda=0$, then $\theta_1=y_1$ annihilates $I_a\ne0$. This contradicts the regularity of $\theta_1$. Hence $x_1$ is a nonzerodivisor. Since $x_1Y=0$, we have $Y=0$. Thus $[II^{-1}]_1=X$, as required.
\end{proof}

For a sequence $\Theta=(\theta_1,\ldots,\theta_n)$ and an $R$-module $N$, set $\Theta N=(\theta_1,\ldots,\theta_n)N$. A sequence $\Theta$ of $d$ linear forms is a \emph{linear system of parameters} if $R/\Theta R$ is Artinian. The algebra $R/\Theta R$ is an \emph{Artinian reduction} of $R$. We allow a field extension when choosing the parameters.

Fix a $\kk$-basis $x_1,\ldots,x_e$ of $R_1$. Let $F=\kk(t_{ij}:1\le i\le d,\ 1\le j\le e)$, where the $t_{ij}$ are algebraically independent over $\kk$, and set $\theta_i=\sum_{j=1}^e t_{ij}x_j$ and $R_F=R\otimes_{\kk}F$. Following Papadakis and Petrotou \cite[Definition~2.2]{PP20}, we call $A_{\mathrm{gen}}=R_F/(\theta_1,\ldots,\theta_d)R_F$ the \emph{generic Artinian reduction} of $R$. The forms $\theta_1,\ldots,\theta_d$ are a linear system of parameters, since such systems form a nonempty open subset of the coefficient space over $\overline{\kk}$.

A standard graded Artinian algebra $A$ has the \emph{strong Lefschetz property} if there is a linear form $\ell\in A_1$ such that, for every pair of integers $0\le a<b$, the map $\times\ell^{b-a}\colon A_a\to A_b$ has maximal rank. Such a form is called a \emph{strong Lefschetz element}.

By \cite[Corollary~4.4.6~(a)]{BH98}, $\indeg\omega_R=d-s$ and $\dim_{\kk}[\omega_R]_{d-s}=h_s$.

To compare entries of $h(R)$, we reduce $R$ and $\omega_R$ modulo the same parameters. The Hilbert function of the first reduction gives $h_0,\ldots,h_s$ starting in degree zero, and that of the second gives $h_s,\ldots,h_0$ starting in degree $d-s$.

\begin{lemma}\label{lem:reduction-pairing}
Let $K/\kk$ be any field extension and let $\Theta$ be any linear system of parameters of $R_{K}=R\otimes_{\kk}K$. Set $A=R_{K}/\Theta R_{K}$ and $M=(\omega_R\otimes_{\kk}K)/\Theta (\omega_R\otimes_{\kk}K)$. For each integer $i$, multiplication induces a perfect pairing $A_i\times M_{d-i}\longrightarrow M_d\cong K$, that is, a $K$-linear isomorphism $M_{d-i}\to\Hom_K(A_i,M_d)$ given by $m\mapsto(a\mapsto am)$. The highest nonzero degree of $M$ is $d$, and $\dim_{K}A_i=\dim_{K}M_{d-i}=h_i$ for every $i$ with $0\le i\le s$. Every nonzero class in $A_i$ pairs nontrivially with the image of some element of $[\omega_R]_{d-i}$. Likewise, every nonzero class in $M_{d-i}$ pairs nontrivially with the image of some element of $R_i$.
\end{lemma}
\begin{proof}
The sequence $\Theta$ is regular because $R_{K}$ is Cohen--Macaulay. Consequently $\Hilb_A(t)=h_R(t)$ and $M(d)\cong\omega_A\cong\Hom_{K}(A,K)$ as graded $A$-modules \cite[Corollary~3.6.14, Theorem~3.6.19~(a), and Exercise~3.3.31]{BH98}. The component of $\Hom_K(A,K)$ in degree $-i$ is $\Hom_K(A_i,K)$, and $A$ acts by $(b\varphi)(a)=\varphi(ba)$. Thus $M_{d-i}$ is the dual of $A_i$, and the pairing is evaluation since $(b\varphi)(1)=\varphi(b)$. Since $A_0=K$ and $A_i=0$ for every $i<0$, the highest nonzero degree of $M$ is $d$. Finally, the images of $\kk$-bases of $[\omega_R]_{d-i}$ and $R_i$ span $M_{d-i}$ and $A_i$, respectively, over $K$. For each nonzero class, at least one such basis element therefore has nonzero product with it.
\end{proof}

Under $(\natur)$ we may choose the parameters in $[\tr_R(\omega_R)]_1$ after extending $\kk$ to $\kk(t)$ if necessary \cite[Remark~3.3~(1), (3), (4)]{M26}. The highest nonzero degrees of the two reductions in Lemma~\ref{lem:reduction-pairing} then give the following inclusions.

\begin{lemma}\label{lem:trace-degree-bounds}
If $R$ satisfies $(\natur)$, then $R_{>s}\subseteq[\tr_R(\omega_R)]_1R$ and $[\omega_R]_{>d}\subseteq[\tr_R(\omega_R)]_1\omega_R$.
\end{lemma}
\begin{proof}
Choose such a system of parameters $\Theta$. By Lemma~\ref{lem:reduction-pairing}, $R/\Theta R$ and $\omega_R/\Theta\omega_R$ vanish in degrees greater than $s$ and $d$, respectively. Each $\theta_i$ belongs to $[\tr_R(\omega_R)]_1$, so the inclusions follow. They descend to $\kk$ by faithful flatness.
\end{proof}

The following lemma gives the two inclusions used in the proof of Theorem~\ref{thm:A}.

\begin{lemma}\label{lem:canonical-squares}
Let $I\subsetneq R$ be a graded ideal containing a homogeneous nonzerodivisor and isomorphic to $\omega_R$ as a graded module. Set $r=\indeg I$, so that $s=d-r$. Suppose that at least one of the following conditions holds.
\begin{enumerate}
\item $R$ satisfies $(\natur)$;
\item $R$ is level and anticanonical level, and $[\tr_R(\omega_R)]_1\ne0$;
\item $R$ is level and $r=1$.
\end{enumerate}
Then
\[
\begin{aligned}
 u^2&\in I_rR &&\text{for every }i\text{ with }0\leq 2i\leq s\text{ and every }u\in I_{d-i},\\
 u^2&\in [I^{-1}]_{1-r}R &&\text{for every }i\text{ with }0\leq 2i<s\text{ and every }u\in R_{s-i}.
\end{aligned}
\]
\end{lemma}
\begin{proof}
Since $I$ is proper, $r\ge1$. If $s=0$, then $d=r$ and only the first inclusion occurs, with $u\in I_r$, so $u^2\in I_rR$. We may therefore assume $s\ge1$, in which case $d=s+r\ge2$.

If $R$ satisfies $(\natur)$, then $\sqrt{[\tr_R(I)]_1R}=\m_R$, so $[\tr_R(I)]_1R$ has height $d\geq2$. By Lemmas~\ref{lem:linear-trace-product} and~\ref{lem:trace-degree-bounds}, we have $[\tr_R(I)]_1=I_r[I^{-1}]_{1-r}$, $I_{>d}\subseteq[\tr_R(I)]_1I$, and $R_{>s}\subseteq[\tr_R(I)]_1R$. If $0\leq2i\leq s$, then every $u\in I_{d-i}$ has square in $I_{2(d-i)}$, where $2(d-i)\geq 2d-s=d+r>d$. Hence $u^2\in[\tr_R(I)]_1I\subseteq I_rR$, since $I^{-1}I\subseteq R$. If $0\leq2i<s$, then $2(s-i)\geq s+1$. Thus every $u\in R_{s-i}$ satisfies $u^2\in[\tr_R(I)]_1R\subseteq[I^{-1}]_{1-r}R$, since $I_r\subseteq R$.

In case~\textup{(2)}, a Gorenstein ring satisfies $(\natur)$, so we may assume that $II^{-1}$ is proper. The generating degrees of $I$ and $I^{-1}$ sum to one, since their product is generated in a single degree and $[II^{-1}]_1\ne0$. Hence $I=I_rR$ and $I^{-1}=[I^{-1}]_{1-r}R$. The first inclusion follows because $u^2\in I$ for every $u\in I$. The second follows from $R\subseteq I^{-1}$.

Finally, under condition~\textup{(3)}, we have $r=1$, $I=I_1R$, and $1\in[I^{-1}]_0$. Hence $u^2\in I_1R$ for every $u\in I$, while $R\subseteq[I^{-1}]_0R$.
\end{proof}

\section{Frobenius splittings and Ehrhart rings}\label{sec:frobenius}

In this section, we prove the alternating inequalities in Theorem~\ref{thm:A} by comparing the reductions of an Ehrhart ring and its canonical ideal modulo the same linear parameters. We use a Frobenius splitting in characteristic two to establish the required injectivity for two multiplication maps between these reductions. The same maps give the Lefschetz assertion in Corollary~\ref{cor:ehrhart-consequences}.

In characteristic two, the \emph{Frobenius endomorphism} of a ring $R$ is the map $a\mapsto a^2$. A \emph{Frobenius splitting} is an additive map $\sigma:R\to R$ such that $\sigma(a^2b)=a\sigma(b)$ for all $a,b\in R$ and $\sigma(1)=1$. In particular, $\sigma(a^2)=a$. An ideal $I$ is \emph{compatibly split} if $\sigma(I)\subseteq I$ \cite[Definition~1.1.3]{BK05}. If $a^{2^e}\in I$ for some $e\geq0$, then repeated application of $\sigma$ gives $a\in I$. Thus every compatibly split ideal is radical, and taking $I=0$ shows that $R$ is reduced.

\begin{example}\label{ex:ehrhart-splitting}
Let $P\subset\RR^q$ be a lattice polytope. Set $C=\RR_{\geq0}(P\times\{1\})$ and $L=\operatorname{span}_{\RR}(C)\cap\ZZ^{q+1}$. The Ehrhart ring $A_P=\kk[C\cap L]$ is a normal affine monoid ring, so it is $F$-pure in positive characteristic \cite[Proposition~6.2 and Corollary~6.3]{BLR06}. Over $\kk=\mathbb F_2$, the standard toric Frobenius splitting \cite[Proposition~3.1]{P09} is given on monomials by
\[
 \sigma_P(X^{\boldsymbol a})=
 \begin{cases}X^{\boldsymbol a/2},&\boldsymbol a\in2L,\\0,&\boldsymbol a\notin2L.\end{cases}
\]
If $\boldsymbol a\in C\cap2L$, then $\boldsymbol a/2\in C\cap L$. Since $\sigma_P(1)=1$ and $\sigma_P(X^{2\boldsymbol b+\boldsymbol a})=X^{\boldsymbol b}\sigma_P(X^{\boldsymbol a})$, the $\mathbb F_2$-linear extension of $\sigma_P$ is a Frobenius splitting. The monomials $X^{\boldsymbol a}$ with $\boldsymbol a\in\relint(C)\cap L$ span a proper ideal $I$. By \cite[Theorem~6.3.5]{BH98}, $I\cong\omega_{A_P}$ as graded $A_P$-modules for the Ehrhart grading. The splitting preserves $I$, since $\boldsymbol a/2\in\relint(C)\cap L$ whenever $\boldsymbol a\in\relint(C)\cap2L$ \cite[proof of Proposition~3.2]{P09}.
\end{example}

\begin{example}\label{ex:stanley-reisner-splitting}
For a simplicial complex $\Delta$ on the vertices $1,\ldots,n$, set $S=\mathbb F_2[x_1,\ldots,x_n]$. On the monomials of $S$, define
\[
\sigma(x^{\boldsymbol a})=\begin{cases}
x^{\boldsymbol a/2}, & \text{if every exponent is even},\\
0, & \text{if some exponent is odd}.
\end{cases}
\]
Multiplying by the square of a monomial adds even exponents, so the $\mathbb F_2$-linear extension satisfies $\sigma(u^2m)=u\sigma(m)$ for all $u,m\in S$. Since $\sigma(1)=1$, this is a Frobenius splitting of $S$. For every monomial $m$ with $\sigma(m)\ne0$, the same variables divide $m$ and $\sigma(m)$. Thus $\sigma$ preserves every ideal generated by squarefree monomials, including the Stanley--Reisner ideal of $\Delta$. The induced Frobenius splitting of $\mathbb F_2[\Delta]$ preserves every ideal generated by squarefree monomials for the same reason.
\end{example}

The following is the key lemma in the proof of Theorem~\ref{thm:A}.

\begin{lemma}\label{lem:frobenius}
Let $R$ be a standard graded algebra over a field $\kk$ of characteristic two. Let $\sigma$ be a Frobenius splitting of $R$, and let $\mathfrak a$ be an ideal with $\sigma(\mathfrak a)\subseteq\mathfrak a$. Let $U\subseteq Q^*(R)$ be a nonzero finite-dimensional $\kk$-subspace consisting of homogeneous elements of degree $\delta$. Choose bases $x_1,\ldots,x_e$ of $R_1$ and $g_1,\ldots,g_b$ of $U$, and fix $n\geq1$. Let $T$ denote a tuple of algebraically independent variables over $\kk$ containing the entries of an $n\times e$ matrix $(t_{ij})$ and the variables $z_1,\ldots,z_b$. Over $K=\kk(T)$, set $\theta_i=\sum_{j=1}^{e}t_{ij}x_j$, $q=\sum_{\ell=1}^{b}z_\ell g_\ell$, and $\Theta=(\theta_1,\ldots,\theta_n)$. Here $R_K=R\otimes_{\kk}K$ and $\mathfrak a_K=\mathfrak a\otimes_{\kk}K$.

Let $f\in R_K$, and let $u$ be any nonzero homogeneous element of $Q^*(R)$. If $u^2\in UR$ and $2\deg u-\delta\geq n$, then $uf\in\Theta \mathfrak a_K$ whenever $qf^2\in\Theta \mathfrak a_K$.
\end{lemma}
\begin{proof}
By hypothesis, $qf^2=\sum_{i=1}^n\theta_i f_i$ for some $f_1,\ldots,f_n\in\mathfrak a_K$. Multiplying by $\theta_1\cdots\theta_n$ gives
\begin{equation}\label{eq:derivative-squares}
 q\theta_1\cdots\theta_n f^2=\sum_{i=1}^n\theta_i^2 F_i,\qquad F_i=\Bigl(\prod_{j\ne i}\theta_j\Bigr)f_i\in\mathfrak a_K.
\end{equation}
For each variable $v$ in $T$, let $\partial_v$ denote formal differentiation with respect to $v$ on $Q^*(R)\otimes_{\kk}K$. It sends $v$ to $1$ and vanishes on $Q^*(R)$ and on the other variables in $T$. The product rule determines $\partial_v$ on rational functions as well, since $\partial_v(h^{-1})=-h^{-2}\partial_v h$ for every nonzero $h\in\kk[T]$. If $c_\nu\in K$ and $a_\nu\in\mathfrak a$, then $\partial_v(\sum_\nu c_\nu a_\nu)=\sum_\nu(\partial_v c_\nu)a_\nu\in\mathfrak a_K$, since $\partial_v c_\nu\in K$. In characteristic two, $\partial_v(f^2)=2f\,\partial_v f=0$, even if $f$ depends on $v$.

Fix $1\leq\ell\leq b$ and $j_1,\ldots,j_n\in\{1,\ldots,e\}$. Apply $D=\partial_{z_\ell}\partial_{t_{1j_1}}\cdots\partial_{t_{nj_n}}$ to~\eqref{eq:derivative-squares}. On the left, the derivatives annihilate $f^2$. Among the remaining factors, $z_\ell$ occurs only in $q$, with $\partial_{z_\ell}q=g_\ell$, and $t_{ij_i}$ occurs only in $\theta_i$, with $\partial_{t_{ij_i}}\theta_i=x_{j_i}$.

On the right, each coefficient derivative $\partial$ satisfies $\partial(\theta_i^2)=2\theta_i\partial(\theta_i)=0$ in characteristic two. The product rule therefore gives $\partial(\theta_i^2 F)=\theta_i^2\partial F$ for every $F$. Applying this at each step yields $D(\theta_i^2 F_i)=\theta_i^2D(F_i)$, so~\eqref{eq:derivative-squares} becomes
\begin{equation}\label{eq:derivative-products}
 g_\ell x_{j_1}\cdots x_{j_n}f^2=\sum_{i=1}^n\theta_i^2D(F_i)\in(\theta_1^2,\ldots,\theta_n^2)\mathfrak a_K.
\end{equation}

By homogeneity, $u^2\in UR$ gives $u^2=\sum_{\ell=1}^b g_\ell w_\ell$ with $w_\ell\in R_{2\deg u-\delta}$. Since $R$ is standard graded and $2\deg u-\delta\geq n$, each $w_\ell$ is an $R$-linear combination of the products $x_{j_1}\cdots x_{j_n}$. Taking the corresponding linear combinations in~\eqref{eq:derivative-products} gives
\[
 (uf)^2=\sum_{\ell=1}^b g_\ell w_\ell f^2\in(\theta_1^2,\ldots,\theta_n^2)\mathfrak a_K.
\]

We now use the splitting to pass from $(uf)^2$ to $uf$. On $R[T]$, set
\[
 \widetilde\sigma\!\left(\sum_{\nu}a_{\nu}T^{\nu}\right)=\sum_{\mu}\sigma(a_{2\mu})T^\mu.
\]
Thus we keep the terms whose exponents are all even and halve their exponents. Multiplication by $c^2T^{2\alpha}$ gives the factor $cT^\alpha$ after applying $\widetilde\sigma$. Since squares add in characteristic two and $\widetilde\sigma(1)=1$, this defines a splitting. It preserves $\mathfrak a[T]$ because $\sigma(\mathfrak a)\subseteq\mathfrak a$. The localization formula $\widetilde\sigma(h/s)=\widetilde\sigma(hs)/s$ extends the splitting to $R_K$. This extension preserves $\mathfrak a_K$. Localizing further at the homogeneous nonzerodivisors of $R$ gives a splitting of $Q^*(R)\otimes_{\kk}K$ that agrees with it on $R_K$. Choose $b_i\in\mathfrak a_K$ such that $(uf)^2=\sum_i\theta_i^2 b_i$. Applying $\widetilde\sigma$ gives
\[
 uf=\widetilde\sigma\bigl((uf)^2\bigr)=\sum_{i=1}^n\theta_i\widetilde\sigma(b_i)\in\Theta \mathfrak a_K.\qedhere
\]
\end{proof}

\begin{remark}
The proof of Lemma~\ref{lem:frobenius} was inspired by the differentiation of the volume map in \cite[Lemma~5.1]{APP25}. Here we differentiate an equation expressing $qf^2\in\Theta\mathfrak a_K$. Since the derivatives act trivially on $Q^*(R)$, the calculation also applies when the $g_\ell$ lie in a fractional ideal.
\end{remark}

Let $I$ be a graded ideal containing a nonzerodivisor. After reduction modulo the parameters, multiplication by an element of $I$ gives a map from the reduction of $R$ to that of $I$. Multiplication by an element of $I^{-1}$ gives a map in the opposite direction. We use Lemma~\ref{lem:frobenius} to prove injectivity for the two maps below.

\begin{proposition}\label{prop:split-injections}
Let $R$ be a standard graded algebra over a field $\kk$ of characteristic two, with a Frobenius splitting $\sigma$. Let $I$ be a graded ideal containing a homogeneous nonzerodivisor, and suppose that $\sigma(I)\subseteq I$. Set $d=\dim R$, $r=\indeg I$, and $s=d-r$. Choose $\kk$-bases $x_1,\ldots,x_e$ of $R_1$, $g_1,\ldots,g_b$ of $I_r$, and $v_1,\ldots,v_c$ of $[I^{-1}]_{1-r}$. Let $F=\kk(t_{j\ell})$ and $K=F(z_1,\ldots,z_b,w_1,\ldots,w_c)$, where all the coefficients are algebraically independent over $\kk$. For each $j$ with $1\leq j\leq d$, set $\theta_j=\sum_{\ell=1}^e t_{j\ell}x_\ell$. Also set
\[
 q=\sum_{\ell=1}^b z_\ell g_\ell,\qquad q'=\sum_{\ell=1}^c w_\ell v_\ell.
\]
If $[I^{-1}]_{1-r}=0$, take $q'=0$. Set $R_K=R\otimes_{\kk}K$ and $I_K=I\otimes_{\kk}K$. The forms $\Theta=(\theta_1,\ldots,\theta_d)$ are a linear system of parameters of $R_K$. Set $A=R_K/\Theta R_K\cong A_{\mathrm{gen}}\otimes_F K$ and $M=I_K/\Theta I_K$. In each assertion $f$ is any homogeneous representative of the class $\bar f$ over $K$, and $u$ is chosen over $\kk$.
\begin{enumerate}
\item Let $0\leq2i\leq s$. Suppose that for every $0\neq\bar f\in A_i$ there is an element $u\in I_{d-i}$ such that $u\bar f\neq0$ in $M_d$ and $u^2\in I_r R$. Then $qf^2\notin\Theta I_K$ for every such $\bar f$. In particular, multiplication by $q$ gives an injection $A_i\longrightarrow M_{r+i}$.
\item Let $0\leq2i<s$. Suppose that for every $0\neq\bar f\in M_{r+i}$ there is an element $u\in R_{s-i}$ such that $u\bar f\neq0$ in $M_d$ and $u^2\in[I^{-1}]_{1-r}R$. Then $q'f^2\notin\Theta I_K$ for every such $\bar f$. In particular, multiplication by $q'$ gives an injection $M_{r+i}\longrightarrow A_{i+1}$.
\end{enumerate}
In particular, if both hypotheses hold at $i$, then $\dim_K A_i\leq\dim_K M_{r+i}\leq\dim_K A_{i+1}$.
\end{proposition}
\begin{proof}
If $d=0$, then $R=I=\kk$, and the assertions follow from $qf^2\ne0$ for nonzero scalars $q$ and $f$. Suppose $d\geq1$. The $d$-tuples of linear parameters form a nonempty open set over $\overline{\kk}$. Since the $t_{j\ell}$ are algebraically independent, $\Theta$ is a system of parameters. The two multiplication maps are well defined because $qR_K\subseteq I_K$ and $q'I_K\subseteq R_K$.

For \textup{(1)}, let $0\neq\bar f\in A_i$ and choose a representative $f\in[R_K]_i$. Suppose $qf^2\in\Theta I_K$ and choose $u$ as in \textup{(1)}. Since $2i\leq s$, we have $2\deg u-r=2(d-i)-r=d+s-2i\geq d$. Together with $u^2\in I_rR$, this allows us to apply Lemma~\ref{lem:frobenius} with $\mathfrak a=I$ and $U=I_r$. It follows that $uf\in\Theta I_K$. Thus $u\bar f=0$ in $M_d$, contrary to the choice of $u$.

For \textup{(2)}, let $0\neq\bar f\in M_{r+i}$ and choose a representative $f\in[I_K]_{r+i}$. Suppose $q'f^2\in\Theta I_K$ and choose $u$ as in \textup{(2)}. Its nonzero square belongs to $[I^{-1}]_{1-r}R$, so $[I^{-1}]_{1-r}\neq0$. Since $2i\leq s-1$, we have $2\deg u-(1-r)=2(s-i)-1+r=d+s-1-2i\geq d$. Apply Lemma~\ref{lem:frobenius} with $\mathfrak a=I$, $U=[I^{-1}]_{1-r}$, and $q'$ in place of $q$. We again get $u f\in\Theta I_K$, contrary to the choice of $u$.

Finally, if $q\bar f=0$ in~\textup{(1)}, multiply $qf\in\Theta I_K$ by $f\in R_K$. If $q'\bar f=0$ in~\textup{(2)}, multiply $q'f\in\Theta R_K$ by $f\in I_K$. Both contradict the respective square nonvanishing.

More generally, the proof of \textup{(1)} applies in any degree if, for each nonzero source class, one can choose $u$ with nonzero product in $M$, $u^2\in I_rR$, and $2\deg u-r\ge d$. The proof of \textup{(2)} applies under the analogous assumptions $u^2\in[I^{-1}]_{1-r}R$ and $2\deg u-(1-r)\ge d$.
\end{proof}

For a Cohen--Macaulay algebra, the pairing in Lemma~\ref{lem:reduction-pairing} gives the nonzero products required in Proposition~\ref{prop:split-injections}. We obtain the following criterion.

\begin{theorem}\label{thm:split-canonical}
Let $R$ be a standard graded Cohen--Macaulay algebra of dimension $d$ over a field $\kk$ of characteristic two. Let $I\subseteq R$ be a graded ideal isomorphic to $\omega_R$ as a graded module, and suppose that $R$ admits a Frobenius splitting $\sigma$ with $\sigma(I)\subseteq I$. Let $r=\indeg I$ and $h(R)=(h_0,\ldots,h_s)$, so that $s=d-r$.

For a fixed integer $i$, consider the conditions
\begin{alignat}{2}
 u^2&\in I_rR &&\qquad\text{for every }u\in I_{d-i},\label{eq:canonical-square-membership}\\
 u^2&\in [I^{-1}]_{1-r}R &&\qquad\text{for every }u\in R_{s-i}.\label{eq:inverse-square-membership}
\end{alignat}
If $0\leq 2i\leq s$ and \eqref{eq:canonical-square-membership} holds, then $h_i\leq h_{s-i}$. If $0\leq 2i<s$ and \eqref{eq:inverse-square-membership} holds, then $h_{s-i}\leq h_{i+1}$. If both conditions hold for every $i$ in their respective ranges, then $h(R)$ is alternatingly increasing, and hence unimodal.
\end{theorem}
\begin{proof}
Since $R$ is reduced, its canonical ideal $I$ contains a homogeneous nonzerodivisor. Take $A$ and $M$ as in Proposition~\ref{prop:split-injections}. Lemma~\ref{lem:reduction-pairing} gives, for every nonzero class in $A_i$, an element $u\in I_{d-i}$ pairing nontrivially with it. Likewise, each nonzero class in $M_{r+i}$ pairs nontrivially with some $u\in R_{s-i}$. Conditions~\eqref{eq:canonical-square-membership} and~\eqref{eq:inverse-square-membership}, respectively, now give the two injections of Proposition~\ref{prop:split-injections}. Their dimensions give the desired inequalities, since $\dim_K A_j=h_j$ and $\dim_K M_{r+i}=h_{s-i}$.
\end{proof}

We are now ready to prove Theorem~\ref{thm:A}. The lattice points in the cone over $P$ determine which monomials belong to the canonical trace ideal of $A_P$, so the conditions on this ideal are independent of the characteristic. Since the remaining hypotheses and the $h$-vector also depend only on $P$, we may apply Theorem~\ref{thm:split-canonical} over $\mathbb F_2$.

\begin{proof}[Proof of Theorem~\ref{thm:A}]
Let $I$ denote the ideal of interior monomials in $A_P$. The exponent sets of $I$ and $I^{-1}$ depend only on $P$ \cite[Proposition~4.1 and Corollary~4.2]{HMP19}. Thus the degrees of their minimal monomial generators, and hence levelness and anticanonical levelness, are independent of the field. The condition $[II^{-1}]_1\neq0$ means that $II^{-1}$ contains a degree-one monomial, and $(\natur)$ means that every degree-one monomial of $A_P$ has a positive power in $[II^{-1}]_1A_P$. Both conditions depend only on these exponent sets. We may thus take $\kk=\mathbb F_2$, where the monomial splitting $\sigma_P$ preserves $I$ by Example~\ref{ex:ehrhart-splitting}. In case~\textup{(3)}, an interior lattice point of $P$ gives $\indeg I=1$. Thus in each case Lemma~\ref{lem:canonical-squares} gives \eqref{eq:canonical-square-membership} and \eqref{eq:inverse-square-membership} in their respective ranges. Theorem~\ref{thm:split-canonical} gives the alternating inequalities over $\mathbb F_2$. Since the Hilbert series counts lattice points in the dilates of $P$, the same inequalities hold over the original field $\kk$.
\end{proof}

\begin{remark}\label{rem:hibi-h-vector}
Let $Q$ be a finite poset and let $L=J(Q)$ be its lattice of order ideals. The Hibi ring of $L$ is the standard graded Ehrhart ring of the order polytope of $Q$. Its $h$-polynomial agrees with that of the order complex of $L$, since both the Hibi ring and the Stanley--Reisner ring of the order complex have bases indexed in each degree by multichains in $L$. Suppose that each connected component of $Q$ is pure and that the ranks of any two components differ by at most one. Then the Hibi ring is nearly Gorenstein \cite[Theorem~5.4]{HHS19}. Hence the $h$-vector of $L$ is alternatingly increasing by Theorem~\ref{thm:A}~\textup{(1)}.
\end{remark}

\begin{remark}
By \cite[Proposition~4.1]{M26}, the Ehrhart ring of an IDP lattice polytope $P$ satisfies $(\natur)$ exactly when its canonical trace ideal contains $\mathbf{x}^w t$ for every vertex $w$ of $P$.

We call $R$ \emph{locally Gorenstein} if $R_{\mathfrak p}$ is Gorenstein for every prime $\mathfrak p\ne\m_R$. Every ring satisfying $(\natur)$ is locally Gorenstein because the canonical trace ideal defines the non-Gorenstein locus \cite[Lemma~2.1]{HHS19}. Even a locally Gorenstein Ehrhart ring can have a non-unimodal $h$-vector. Ferroni constructed a smooth IDP polytope $P$ with non-unimodal $h^*(P)$ \cite[Theorem~4.1]{F26}. Since $A_P$ is standard graded and $\operatorname{Proj} A_P$ is smooth, $A_P$ is regular on the punctured spectrum.
\end{remark}

Under the hypotheses for the alternating inequalities in Theorem~\ref{thm:split-canonical}, one can choose $q$ and $q'$ so that $\overline{qq'}$ is a strong Lefschetz element of the generic Artinian reduction.

\begin{theorem}\label{thm:split-strong-lefschetz}
Let $R$ be a $d$-dimensional standard graded Cohen--Macaulay algebra over a field $\kk$ of characteristic two. Let $I\subseteq R$ be a graded ideal isomorphic to $\omega_R$ as a graded $R$-module. Suppose that $R$ has a Frobenius splitting $\sigma$ such that $\sigma(I)\subseteq I$. Set $r=\indeg I$ and $s=d-r$. Assume that
\[
\begin{aligned}
u^2&\in I_rR &&\text{for every $i$ with $0\leq 2i\leq s$ and every $u\in I_{d-i}$},\\
u^2&\in[I^{-1}]_{1-r}R &&\text{for every $i$ with $0\leq 2i<s$ and every $u\in R_{s-i}$}.
\end{aligned}
\]
Let $A_{\mathrm{gen}}$ be the generic Artinian reduction of $R$ over $F=\kk(t_{ij})$, where the $t_{ij}$ are the coefficients of the generic linear system of parameters. In fixed $\kk$-bases of $I_r$ and $[I^{-1}]_{1-r}$, respectively, take linear combinations $q$ and $q'$ whose coefficients are jointly algebraically independent over $F$. Let $K$ be the field obtained by adjoining these coefficients to $F$. Then the image of $qq'$ is a strong Lefschetz element of $A=A_{\mathrm{gen}}\otimes_F K$. Moreover, one can choose the coefficients of $q$ and $q'$ in $F$ so that the image of their product is a strong Lefschetz element of $A_{\mathrm{gen}}$.
\end{theorem}
\begin{proof}
If $s=0$, then $h_R(t)=1$, so both $A_{\mathrm{gen}}$ and $A$ are concentrated in degree zero. For every choice of coefficients, the degree-one element $\overline{qq'}$ is therefore zero. It is a strong Lefschetz element since every map in the definition has zero target and hence maximal rank. Now assume $s\ge1$. Let $\Theta$ be the generic linear system of parameters defining $A_{\mathrm{gen}}$, and set $I_K=I\otimes_{\kk}K$ and $M=I_K/\Theta I_K$. Consider the sequence of multiplication maps
\[
 A_0\xrightarrow{\times q}M_r\xrightarrow{\times q'}A_1\longrightarrow\cdots\longrightarrow A_s\xrightarrow{\times q}M_d.
\]
Set $D_{2i}=A_i$ and $D_{2i+1}=M_{r+i}$ for each $i$ with $0\le i\le s$. Alternating multiplication by $q$ and $q'$ defines a degree-one linear map $\varphi$ on $\bigoplus_{j=0}^{2s+1}D_j$, with $\varphi(D_{2s+1})=0$. For $0\le j\le 2s+1$, the perfect pairing $D_j\times D_{2s+1-j}\to M_d$ from Lemma~\ref{lem:reduction-pairing} is given by $\langle x,y\rangle=xy$. For every $0\le j\le 2s$ and all $x\in D_j$, $y\in D_{2s-j}$, we have $\langle\varphi x,y\rangle=\langle x,\varphi y\rangle$. If $j$ is even, $x,y\in A$ and both sides equal $\bar q(xy)$, where $\bar q$ is the class of $q$ in $M_r$. If $j$ is odd, take representatives $\tilde x,\tilde y\in I_K$ of $x,y\in M$. The two sides are $(q'x)y$ and $x(q'y)$, respectively. Since $q'I_K\subseteq R_K$, both are the class of $(q'\tilde x)\tilde y=\tilde x(q'\tilde y)$ modulo $\Theta I_K$.

For $0\le 2i\le s$, Lemma~\ref{lem:reduction-pairing} gives, for every nonzero $x\in A_i$, an element $u\in I_{d-i}$ such that $\bar u x\ne0$ in $M_d$. For $0\le 2i<s$, the same lemma gives, for every nonzero $y\in M_{r+i}$, an element $v\in R_{s-i}$ such that $\bar v y\ne0$ in $M_d$, since $d-(r+i)=s-i$. Here $\bar u$ and $\bar v$ are the images in $M$ and $A$, respectively. The hypotheses on squares give $u^2\in I_rR$ and $v^2\in[I^{-1}]_{1-r}R$ in the respective ranges. Proposition~\ref{prop:split-injections}\,\textup{(1)} and \textup{(2)} therefore give the injectivity of $\varphi\colon D_j\to D_{j+1}$ for every $0\le j\le s$.

The middle map $D_s\to D_{s+1}$ is $A_i\xrightarrow{\times q}M_{r+i}$ when $s=2i$, and $M_{r+i}\xrightarrow{\times q'}A_{i+1}$ when $s=2i+1$. Take any nonzero $z\in D_s$ and choose a homogeneous representative $f$. The element $\langle z,\varphi z\rangle\in M_d$ is represented by $qf^2$ if $s$ is even and by $q'f^2$ if $s$ is odd. Proposition~\ref{prop:split-injections}\,\textup{(1)} gives $qf^2\notin\Theta I_K$ when $s$ is even, and Proposition~\ref{prop:split-injections}\,\textup{(2)} gives $q' f^2\notin\Theta I_K$ when $s$ is odd. Thus $\langle z,\varphi z\rangle\ne0$.

Fix $j$ with $0\le j\le s$. For every nonzero $x\in D_j$, we have $z=\varphi^{s-j}x\ne0$ and
\[
 \langle x,\varphi^{2s+1-2j}x\rangle=\langle\varphi^{s-j}x,\varphi^{s+1-j}x\rangle=\langle z,\varphi z\rangle\ne0.
\]
Thus $\varphi^{2s+1-2j}:D_j\to D_{2s+1-j}$ is injective, and hence an isomorphism by duality. Now fix any pair $0\le m\le n\le 2s+1$. If $m+n\le2s+1$, the isomorphism starting in $D_m$ factors through $\varphi^{n-m}\colon D_m\to D_n$, so this map is injective. If $m+n>2s+1$, the isomorphism from $D_{2s+1-n}$ to $D_n$ factors through $\varphi^{n-m}\colon D_m\to D_n$, so this map is surjective.

On $A_a=D_{2a}$, the map $\varphi^{2(b-a)}$ is multiplication by $\ell^{b-a}$, where $\ell=\overline{qq'}$. The map $A_a\to A_b$ is injective if $a+b\le s$ and surjective if $a+b\ge s+1$. Thus it has rank $\min\{h_a,h_b\}$ for every $0\le a<b\le s$, and $\ell$ is a strong Lefschetz element over $K$. In fixed $F$-bases, choose one nonzero maximal minor for each of these maps. Their product is a nonzero polynomial over $F$ in the coefficients of $q$ and $q'$. Since $F$ is infinite, we may choose the coefficients in $F$ so that this product is nonzero. The resulting $\overline{qq'}$ is a strong Lefschetz element of $A_{\mathrm{gen}}$.
\end{proof}

Combining Lemma~\ref{lem:canonical-squares} with Theorems~\ref{thm:split-canonical} and~\ref{thm:split-strong-lefschetz} gives the following corollary.

\begin{corollary}\label{cor:split-applications}
Let $R$ be a standard graded Cohen--Macaulay algebra over a field of characteristic two, with a Frobenius splitting $\sigma$ and a proper graded ideal $I\cong\omega_R$ such that $\sigma(I)\subseteq I$, where the isomorphism preserves degrees. If one of the following conditions holds, then $h(R)$ is alternatingly increasing and the generic Artinian reduction of $R$ has the strong Lefschetz property:
\begin{enumerate}
\item $R$ satisfies $(\natur)$;
\item $R$ is level and anticanonical level, and $[\tr_R(\omega_R)]_1\neq0$;
\item $R$ is level and $\indeg\omega_R=1$.
\end{enumerate}
\end{corollary}
\begin{proof}
Since $R$ is reduced, $I$ contains a homogeneous nonzerodivisor. In each of the three cases, Lemma~\ref{lem:canonical-squares} gives \eqref{eq:canonical-square-membership} and \eqref{eq:inverse-square-membership} in their respective ranges. Theorem~\ref{thm:split-canonical} gives the alternating inequalities, and Theorem~\ref{thm:split-strong-lefschetz} gives the strong Lefschetz property.
\end{proof}

We apply the multiplication maps to the Ehrhart ring of an IDP lattice polytope. Recall that an \emph{$O$-sequence} is the Hilbert function of a standard graded algebra over a field.

\begin{corollary}\label{cor:ehrhart-consequences}
Let $P$ be an IDP lattice polytope with $h^*$-vector $h^*(P)=(h_0,\ldots,h_s)$, and let $R=A_P$ be its Ehrhart ring over a field $\kk$. Set $h_j=0$ for $j<0$ or $j>s$. Suppose that at least one of the following conditions holds.
\begin{enumerate}
\item $R$ satisfies $(\natur)$;
\item $R$ is level and anticanonical level, and $[\tr_R(\omega_R)]_1\ne0$;
\item $R$ is level and $P$ has an interior lattice point.
\end{enumerate}
Then the following assertions hold.
\begin{enumerate}
\renewcommand{\theenumi}{\alph{enumi}}
\renewcommand{\labelenumi}{\textup{(\theenumi)}}
\item\label{item:ehrhart-lefschetz} If $\kk$ has characteristic zero or two, the generic Artinian reduction of $R$ has the strong Lefschetz property.
\item\label{item:ehrhart-osequences} Both of the following sequences are $O$-sequences.
\[
 \bigl(\max\{h_j-h_{j-1},0\}\bigr)_{j\ge0},\qquad
 \bigl(\max\{h_j-h_{s+1-j},0\}\bigr)_{j\ge0}.
\]
\end{enumerate}
\end{corollary}
\begin{proof}
If $s=0$, the generic Artinian reduction is a field, so $\ell=0$ is a strong Lefschetz element and both sequences are $(1,0,\ldots)$. Assume $s\ge1$.

\smallskip
\noindent\textup{(\ref{item:ehrhart-lefschetz})} Let $I$ be the canonical ideal spanned by the interior monomials as in Example~\ref{ex:ehrhart-splitting}, and set $r=\indeg I$. Use independent variables $\mathbf t$ as coefficients of the parameters $\Theta$ in a monomial basis of $R_1$. Set $F=\kk(\mathbf t)$ and let $A=R_F/\Theta R_F$ be the generic Artinian reduction. In monomial bases of $I_r$ and $[I^{-1}]_{1-r}$, respectively, take linear combinations $q$ and $q'$ with independent coefficients $\mathbf z$ and $\mathbf w$. Set $K=F(\mathbf z,\mathbf w)$ and $A_K=A\otimes_F K$. For $\kk=\mathbb F_2$, Lemma~\ref{lem:canonical-squares} and Theorem~\ref{thm:split-strong-lefschetz} give the strong Lefschetz element $\ell=\overline{qq'}$ of $A_K$.

For each pair $0\le a<b\le s$, let $B_b$ be the matrix whose columns are the coordinate vectors of the products $\theta_j m$ in the monomial basis of $R_b$, with $m$ ranging over the monomials of degree $b-1$. Append the columns of multiplication by $(qq')^{b-a}$ from $R_a$ to $R_b$. The increase in rank is the rank of multiplication by $\overline{qq'}^{\,b-a}$ from $A_a\otimes_F K$ to $A_b\otimes_F K$. In the monomial bases, multiplication adds exponents. Thus the same polynomials in $\mathbf t,\mathbf z,\mathbf w$ with integer coefficients give the entries in every characteristic. Also, $\operatorname{rank}B_b=\dim_{\kk}R_b-h_b$ in every characteristic. For each $(a,b)$ the augmented matrix therefore has a minor of size $\operatorname{rank}B_b+\min\{h_a,h_b\}$ whose determinant is nonzero modulo two. This determinant is nonzero over every field of characteristic zero or two. Since $F=\kk(\mathbf t)$ is infinite, we can choose values in $F$ for $\mathbf z$ and $\mathbf w$ at which the product of the chosen determinants is nonzero. With these coefficients, $\overline{qq'}$ is a strong Lefschetz element of $A$ over $F$.

\smallskip
\noindent\textup{(\ref{item:ehrhart-osequences})} As in the proof of Theorem~\ref{thm:A}, the hypotheses and $h^*(P)$ are independent of the ground field. Take $\kk=\mathbb F_2$ and retain $A_K$, $q$, $q'$ and $\ell=\overline{qq'}$ from part~\textup{(\ref{item:ehrhart-lefschetz})}, with the coefficients of $q$ and $q'$ still algebraically independent over $F$. Since $\ell$ is a strong Lefschetz element, $A_K/(\ell)$ has the first Hilbert function.

Set $I_K=I\otimes_{\kk}K$ and $M=I_K/\Theta I_K$. The degree-$j$ component of the homogeneous ideal $q'M$ is the image of multiplication by $q'$ from $M_{r+j-1}$ to $[A_K]_j$. By Lemma~\ref{lem:reduction-pairing}, the source has dimension $h_{s+1-j}$. For $1\le j\le (s+1)/2$, Lemmas~\ref{lem:reduction-pairing} and~\ref{lem:canonical-squares} give the hypotheses of Proposition~\ref{prop:split-injections}\,\textup{(2)} at $i=j-1$, so this map is injective. For $(s+1)/2<j\le s$, Theorem~\ref{thm:A} gives $h_{j-1}\ge h_j$, so multiplication by $\ell$ from $[A_K]_{j-1}$ to $[A_K]_j$ is surjective. Since $q[A_K]_{j-1}\subseteq M_{r+j-1}$, we obtain $[A_K]_j=\ell[A_K]_{j-1}\subseteq q'M_{r+j-1}\subseteq[A_K]_j$. Multiplication by $q'$ therefore has maximal rank in every degree $1\le j\le s$. Since $M_{r-1}=0$, the degree-zero component of $A_K/q'M$ has dimension one. Hence $A_K/q'M$ has the second Hilbert function. Moreover, $\ell\in q'M$, so there is a natural surjection $A_K/(\ell)\to A_K/q'M$.
\end{proof}

\begin{remark}\label{rem:app-scope}
In the proof of Corollary~\ref{cor:ehrhart-consequences}\,\textup{(\ref{item:ehrhart-lefschetz})}, the required minors are polynomials with integer coefficients, so their nonvanishing modulo two gives nonvanishing in characteristic zero. For a related use of reduction modulo $p$ for the degree pairing of a Stanley--Reisner ring, see \cite[proof of Proposition~3.5]{KLS26}.

In case~\textup{(3)} of Theorem~\ref{thm:A}, the alternating inequalities asked for in \cite[Question~2.21]{SVL13} follow from \cite[Corollary~2.2 and Theorem~2.4]{APP25}, since $s=\dim P$. The strong Lefschetz property in Corollary~\ref{cor:ehrhart-consequences}\,\textup{(\ref{item:ehrhart-lefschetz})} follows in this case from \cite[Theorems~2.1 and~6.1]{APP25} with $j=1$ in Theorem~6.1.
\end{remark}

\begin{example}\label{ex:second-difference}
Consider the projective monomial curve with coordinate ring
\[
 R=\QQ[u^{51},u^{50}v,u^{49}v^2,u^{48}v^3,u^{37}v^{14},u^{36}v^{15},v^{51}],
\]
with all seven generators in degree one. One can check with \texttt{Macaulay2}~\cite{M2} that $R$ is Cohen--Macaulay and $h(R)=(1,5,9,13,12,9,2)$. This vector is concave and alternatingly increasing, with $h_s=2$ and $h_1-h_s=3$. The first sequence in Corollary~\ref{cor:ehrhart-consequences}\,\textup{(\ref{item:ehrhart-osequences})} is the $O$-sequence $(1,4,4,4,0,\ldots)$. The second begins $(1,3,0,1)$, so it cannot be an $O$-sequence. Thus $h(R)$ cannot be the $h^*$-vector of an IDP polytope satisfying the hypotheses of Theorem~\ref{thm:A}.
\end{example}

We next apply Corollary~\ref{cor:split-applications}\,\textup{(3)} to the Stanley--Reisner rings of triangulated balls with induced boundary.
For an $s$-dimensional simplicial complex $\Delta$, let $h(\Delta)=(h_0(\Delta),\ldots,h_{s+1}(\Delta))$ denote its \emph{$h$-vector}, where $h_i(\Delta)$ is the coefficient of $t^i$ in the $h$-polynomial of its Stanley--Reisner ring $\kk[\Delta]$.
The alternating inequalities in the next corollary are known \cite[Section~5.1, (17) and Theorem~5.1]{Ath26}. For the Lefschetz assertion, the injections follow from the argument of Adiprasito and Yashfe~\cite{AY21} explained in \cite[Section~5.1]{Ath26}. The surjections follow by applying the sphere Lefschetz theorem \cite[Corollary~9.10]{PP20} after coning off $\partial\Delta$.

\begin{corollary}\label{cor:ball-induced-boundary}
Let $\Delta$ be a triangulation of an $s$-dimensional ball such that $\partial\Delta$ is an induced subcomplex of $\Delta$. Then $h_{s+1}(\Delta)=0$ and the sequence $h_0(\Delta),\ldots,h_s(\Delta)$ is alternatingly increasing. Moreover, for every field $\kk$ of characteristic two, the generic Artinian reduction of $\kk[\Delta]$ has the strong Lefschetz property.
\end{corollary}
\begin{proof}
For a single vertex, $h(\Delta)=(1,0)$ and the Artinian reduction is a field. For $s\ge1$, the ring $R=\mathbb F_2[\Delta]$ is Cohen--Macaulay. By Hochster's theorem \cite[Theorem~II.7.3]{S96}, the monomials $x_F=\prod_{v\in F}x_v$ for the interior faces $F$ of $\Delta$ generate a graded canonical ideal $I$ of $R$. Since $\partial\Delta$ is induced, every face supported on boundary vertices lies in $\partial\Delta$. Thus each interior face contains an interior vertex, and the variables of the interior vertices generate $I$. Hence $R$ is level with $\indeg\omega_R=1$, and $\deg h_R=\dim R-1=s$.

Since $I$ is generated by variables, the Frobenius splitting of $R$ in Example~\ref{ex:stanley-reisner-splitting} preserves $I$. Corollary~\ref{cor:split-applications}\,\textup{(3)} gives both assertions over $\mathbb F_2$. Extending the parameter field from $\mathbb F_2(\mathbf t)$ to $\kk(\mathbf t)$ gives the Lefschetz assertion over every field $\kk$ of characteristic two.
\end{proof}

\section{\texorpdfstring{$F$}{F}-pure rings in characteristic two}\label{sec:fpure}

In this section, we prove Theorem~\ref{thm:B} and apply it to a Rees algebra and Veronese subrings.
For an $R$-module $M$, let $F_*M$ be the same additive group as $M$ with scalar multiplication $a\cdot F_*m=F_*(a^2m)$, where $F_*m$ denotes the element $m$ in this module. With this action, the Frobenius map $F:R\to F_*R$, $a\mapsto F_*(a^2)$, is $R$-linear. For graded modules we set $\deg(F_*m)=\deg(m)/2$ for homogeneous $m$, so $F$ has degree zero.

We call $R$ \emph{$F$-finite} if $F_*R$ is finitely generated as an $R$-module. For example, $F_*\mathbb F_2[x]$ is generated by $F_*1$ and $F_*x$, corresponding to the even and odd powers of $x$. In several variables one separates the exponents by parity in each variable. The images of the resulting generators work in any finitely generated $\mathbb F_2$-algebra.

A ring homomorphism $R\to S$ is \emph{pure} if $N\to N\otimes_R S$, $n\mapsto n\otimes1$, is injective for every $R$-module $N$. We call $R$ \emph{$F$-pure} if its Frobenius homomorphism is pure \cite[Definition~2.1]{DMV26}.

In this notation, a Frobenius splitting is an $R$-linear map $\sigma:F_*R\to R$ such that $\sigma(F_*1)=1$. Here $R$-linearity means $\sigma(F_*(a^2b))=a\sigma(F_*b)$, so taking $b=1$ gives $\sigma\circ F=\mathrm{id}_R$. For an $F$-finite ring, $F$-purity is equivalent to the existence of such a splitting \cite[Corollary~2.4]{MaPo}. If $R$ is graded, the degree-zero component of a splitting is again a splitting, since it still sends $F_*1$ to $1$.

We prepare two lemmas for the proof of Theorem~\ref{thm:B}. The first constructs a graded canonical ideal $I$ with $\sigma(F_*I)\subseteq I$ for a given splitting $\sigma$, as required in Corollary~\ref{cor:split-applications}.

\begin{lemma}\label{lem:fpure-canonical-splitting}
Let $R$ be an $F$-finite standard graded Cohen--Macaulay algebra of characteristic two that is Gorenstein at every prime of height one. For every Frobenius splitting $\sigma:F_*R\to R$ of degree zero, there is a graded ideal $I\subseteq R$ isomorphic to $\omega_R$ by a degree-preserving map such that $\sigma(F_*I)\subseteq I$.
\end{lemma}
\begin{proof}
The splitting makes $R$ reduced. If $\dim R=0$, all positive-degree elements are nilpotent, so $R=\kk$ and we may take $I=R$. Suppose $\dim R>0$. Since $R$ is generically Gorenstein, we may identify $\omega_R$ with a graded fractional ideal $J\subseteq Q^*(R)$ by a degree-preserving isomorphism \cite[proof of Remark~2.13~(1)]{KM25}.

Set $S=F_*R$. The Frobenius map $F:R\to S$ is finite since $R$ is $F$-finite, and it preserves degrees since $\deg(F_*(x^2))=\deg x$ for homogeneous $x$. For this finite graded map, \cite[Proposition~3.6.12~(b)]{BH98} gives $\omega_S\cong\operatorname{Ext}_R^{\,\dim R-\dim S}(S,\omega_R)$. Since $S$ is $R$ with all degrees divided by two, we have $\dim R=\dim S$. Then $\operatorname{Ext}_R^{\,\dim R-\dim S}(S,\omega_R)=\Hom_R(S,\omega_R)$. The same change of grading on $J\cong\omega_R$ makes $F_*J$ a graded canonical module of $S$. We may therefore choose a graded $S$-module isomorphism $F_*J\cong\Hom_R(S,J)$.

Applying $\Hom_R(-,J)$ to $F$ sends an $R$-linear map $h:S\to J$ to $h\circ F:R\to J$. An $R$-linear map $R\to J$ is determined by its value at $1$, which here is $(h\circ F)(1)=h(F_*1)$. Thus the Frobenius trace is the composite
\[
 \tau\colon F_*J\xrightarrow{\cong}\Hom_R(S,J)\xrightarrow{\,h\mapsto h(F_*1)\,}J.
\]
Both arrows preserve degrees, since the chosen isomorphism is graded and $F_*1$ has degree zero. Localizing at the homogeneous nonzerodivisors extends $\tau$ to $F_*Q^*(R)\to Q^*(R)$.

Since $R$ is Cohen--Macaulay, it satisfies $(S_2)$. Its localizations at minimal primes are fields because $R$ is reduced.
Together with the hypothesis, $R_{\mathfrak p}$ is Gorenstein for every prime $\mathfrak p$ of height at most one. We may therefore apply \cite[Lemma~5.15]{MaPo} with $p=2$ and $e=1$. Since $1-p^e=-1$, we obtain an $R$-module isomorphism $F_*(J^{-1})\cong\Hom_R(F_*R,R)$, where $J^{-1}=(R:J)$.

In terms of $\tau$, this isomorphism sends $F_*c$ to the map $\varphi_c(F_*x)=\tau(F_*(cx))$ \cite[Proposition~2.12 and Example~2.14]{CRS23}. If $c$ is homogeneous of degree~$j$, the map $F_*x\mapsto F_*(cx)$ raises degrees by $j/2$. Since $\tau$ has degree zero, $\varphi_c$ has degree $j/2=\deg(F_*c)$. In particular, the degree-zero splitting $\sigma$ corresponds to some $c\in[J^{-1}]_0$.

Set $I=cJ\subseteq R$. Since $\tau(F_*c)=\sigma(F_*1)=1$ on each factor of $Q^*(R)$, the homogeneous element $c$ is a unit in $Q^*(R)$. Multiplication by $c$ therefore gives a degree-preserving isomorphism $J\cong I$. For every $a\in J$, we have $\sigma(F_*(ca))=\tau(F_*(c^2a))=c\tau(F_*a)\in I$.
\end{proof}

The second reduces the proof to a field finitely generated over $\mathbb F_2$.

\begin{lemma}\label{lem:fpure-field-descent}
Let $R$ be an $F$-pure standard graded Cohen--Macaulay algebra over a field $\kk$ of characteristic two. Suppose that $R$ satisfies $(\natur)$. Then there exist a subfield $K\subseteq\kk$ finitely generated over $\mathbb F_2$ and an $F$-finite standard graded $K$-algebra $R'$ such that $R\cong R'\otimes_K\kk$. Moreover, $R'$ is Cohen--Macaulay and $F$-pure, satisfies $(\natur)$, and $h(R)=h(R')$. If the generic Artinian reduction of $R'$ has the strong Lefschetz property, then so does that of $R$.
\end{lemma}
\begin{proof}
Let $K\subseteq\kk$ be the subfield generated over $\mathbb F_2$ by the coefficients of a finite homogeneous presentation of $R$, and let $R'$ be the algebra defined by this presentation over $K$. Then $R\cong R'\otimes_K\kk$, and $R'$ is $F$-finite since it is essentially of finite type over $\mathbb F_2$ \cite[Exercise~5]{MaPo}. Since $R'\to R$ is faithfully flat, $R'$ is Cohen--Macaulay \cite[Theorem~2.1.10]{BH98} and $F$-pure \cite[Proposition~5.13]{HR76}. The canonical trace ideal commutes with field extension \cite[Remark~3.3~(1)]{M26}, so $(\natur)$ descends to $R'$ by faithful flatness. We have $h(R')=h(R)$. Using the same indeterminates $\mathbf t$ for the generic parameters gives $A\cong A'\otimes_{K(\mathbf t)}\kk(\mathbf t)$ for the generic Artinian reductions of $R$ and $R'$, respectively. Thus every strong Lefschetz element of $A'$ remains one in $A$, since field extension preserves ranks.
\end{proof}

\begin{proof}[Proof of Theorem~\ref{thm:B}]
By Lemma~\ref{lem:fpure-field-descent}, it suffices to treat the $F$-finite case. Choose a degree-zero Frobenius splitting $\sigma$ and set $d=\dim R$. By \cite[Proposition~2.4]{HR76}, $[H^d_{\m_R}(R)]_j=0$ for every $j>0$, so graded local duality \cite[Theorem~3.6.19~(a)]{BH98} gives $\indeg\omega_R\geq0$.

If $d\leq1$, then $\deg h_R=d-\indeg\omega_R\leq1$. Hence $h(R)$ is either $(1)$ or $(1,h_1)$ with $h_1\geq1$. In the generic Artinian reduction, multiplication by any nonzero linear form is injective from degree zero to degree one, and the higher-degree parts vanish. If the reduction is a field, the zero form is a strong Lefschetz element.

Assume $d\geq2$. Condition $(\natur)$ makes the canonical trace ideal the unit ideal away from $\m_R$, so $R$ is Gorenstein at every prime of height one \cite[Lemma~2.1]{HHS19}. Applying Lemma~\ref{lem:fpure-canonical-splitting} to $\sigma$ gives a graded ideal $I\cong\omega_R$ with $\sigma(F_*I)\subseteq I$. If $I$ is proper, Corollary~\ref{cor:split-applications}\,\textup{(1)} applies. If $I=R$, then $r=0$, $s=d$, $I_0R=R$, and $[I^{-1}]_1R=\m_R$. Since $s-i>0$ for $0\leq2i<s$, both square conditions in Theorems~\ref{thm:split-canonical} and~\ref{thm:split-strong-lefschetz} hold.
\end{proof}

In the Gorenstein case over $\mathbb F_2$, \cite[Proposition~4.3]{AKOPP26} also gives the Lefschetz assertion via the application of Fedder's criterion announced in \cite{Hak26}.

The following Rees algebra satisfies the hypotheses of Theorem~\ref{thm:B}.

\begin{example}\label{ex:rees-ordinary-cubic}
Let $C$ be the homogeneous coordinate ring of a smooth ordinary plane cubic over $\kk=\mathbb F_2$. For example, homogenizing the Weierstrass equation $y^2+xy=x^3+x$ in \cite[Section~2, (4)]{Sch22} gives $C=\kk[x,y,z]/(y^2z+xyz+x^3+xz^2)$. Set $A=C[x_1,x_2,x_3]$ and $R=A[\m_A T]$, graded by $\deg(aT^j)=\deg_A a$ for homogeneous $a\in\m_A^j$. The algebra $R$ is generated in degree one by $A_1$ and $A_1T$. Since $A$ is a cubic hypersurface in six variables, $R_{\m_R}$ is Cohen--Macaulay and its canonical trace ideal equals its maximal ideal \cite[Example~2.11]{EY26}. The canonical trace ideal is homogeneous and commutes with localization, so $\tr_R(\omega_R)=\m_R$. Thus $R$ is nearly Gorenstein.

The affine cone over the smooth cubic is regular off the vertex, so $C$ is normal by Serre's criterion. Ordinarity gives a Frobenius splitting of $C$ \cite[Example~2.3.1]{KKP21}. Hence $A$ is a standard graded $F$-pure normal domain, and \cite[Proposition~6.1]{Si00} gives $F$-purity of $R$. As $R$ satisfies $(\natur)$, Theorem~\ref{thm:B} gives the strong Lefschetz property for its generic Artinian reduction.
\end{example}

We next apply Theorem~\ref{thm:B} to Veronese subrings. For a graded algebra $S$ and an integer $r\geq1$, let $S^{(r)}=\bigoplus_{j\geq0}S_{rj}$ denote its \emph{$r$th Veronese subring}, with $[S^{(r)}]_j=S_{rj}$.

\begin{corollary}\label{cor:fpure-veronese}
Let $R$ be an $F$-pure standard graded Cohen--Macaulay algebra over a field of characteristic two that satisfies $(\natur)$. For every $r\geq1$, the $h$-vector of $R^{(r)}$ is alternatingly increasing and its generic Artinian reduction has the strong Lefschetz property. In particular, both conclusions hold for every Veronese subring of an $F$-pure standard graded Gorenstein algebra in characteristic two.
\end{corollary}
\begin{proof}
The projection onto the components of degree divisible by $r$ is $R^{(r)}$-linear, so $R^{(r)}\subseteq R$ is a pure extension. Since $F$-purity descends to pure subrings \cite[Proposition~5.13]{HR76}, $R^{(r)}$ is $F$-pure. It is standard graded and Cohen--Macaulay \cite[Section~3.1]{GW78}, and it satisfies $(\natur)$ by \cite[Theorem~3.6~(1)]{M26}. Theorem~\ref{thm:B} gives both conclusions. If $R$ is Gorenstein, then $\tr_R(\omega_R)=R$, so $R$ satisfies $(\natur)$.
\end{proof}

\begin{example}\label{ex:veronese-ordinary-cubic}
Let $C$ be the homogeneous coordinate ring of the ordinary plane cubic in Example~\ref{ex:rees-ordinary-cubic} and set $A=C[w]$ with $\deg w=1$. This is an $F$-pure standard graded hypersurface with $\omega_A\cong A(-1)$. By \cite[Corollary~3.1.3]{GW78}, the lowest nonzero component of $\omega_{A^{(r)}}$ is $[\omega_{A^{(r)}}]_1\cong A_{r-1}$. For $r\geq2$ this component has dimension greater than one, so Corollary~\ref{cor:fpure-veronese} applies to the non-Gorenstein algebras $A^{(r)}$.
\end{example}

\begin{remark}\label{rem:veronese-ehrhart}
Since the hypotheses and the $h$-vector in Theorem~\ref{thm:A} depend only on $P$, we may work over $\mathbb F_2$. Example~\ref{ex:ehrhart-splitting} gives a Frobenius splitting of $A_P$ that preserves the canonical ideal of interior monomials. In particular, $A_P$ is $F$-pure, so case~\textup{(1)} of Theorem~\ref{thm:A} also follows from Theorem~\ref{thm:B}. Corollary~\ref{cor:split-applications}\,\textup{(2)} and~\textup{(3)} give the remaining cases, since in case~\textup{(3)} an interior lattice point gives $\indeg\omega_{A_P}=1$.
\end{remark}

\section{A numerical criterion for unimodality}\label{sec:domains}

In this section we prove Theorem~\ref{thm:C}. Over a suitable field extension $K/\kk$, we first construct a graded surjection from a Veronese subring of $K[x,y]$ onto an Artinian reduction of $R$. Its Hilbert function gives unimodality. To estimate the middle entry of $h(R)$, we then use a canonical ideal of $R$.

By Lemma~\ref{lem:linear-trace-product}, the degree-one part of the canonical trace ideal of a standard graded Cohen--Macaulay domain is spanned by products $vu$, where $u$ belongs to the lowest-degree component of a canonical fractional ideal and $v$ to the complementary degree of its inverse. The following lemma uses these products to construct the required surjection.

\begin{lemma}\label{lem:product-matrix}
Let $R$ be a standard graded Cohen--Macaulay domain over a field $\kk$. Let $E,W\subseteq Q^*(R)$ be nonzero finite-dimensional $\kk$-subspaces, each concentrated in a single degree, such that $EW\subseteq R_1$. Suppose that $\dim_{\kk}E\geq2$ and $\height((EW)R)\geq2$. Set $p=\dim_{\kk}W$, $K=\overline{\kk}$, and $R_K=R\otimes_{\kk}K$. For every sufficiently general subspace $L\subseteq [R_K]_1$ of codimension $p+1$, there is a graded surjection $K[x,y]^{(p)}\to R_K/(L)$.
\end{lemma}
\begin{proof}
Fix $0\ne v_1\in E$ and take $v\in E\setminus\kk v_1$. The fraction $r=v/v_1$ sends $(v_1W)R$ into $R$. If this ideal contained a regular sequence $a,b$, then $b(ra)=a(rb)$ would give $ra\in aR$, hence $r\in R_0=\kk$. Thus the nonzero ideal $(v_1W)R$ has height one, since $R$ is Cohen--Macaulay. Choose a height-one prime $\mathfrak p$ containing it. The height hypothesis on $(EW)R$ gives $v_2\in E$ and $w\in W$ with $v_2w\notin\mathfrak p$. In particular, $v_2w\notin v_1W$, so $\dim_{\kk}R_1\geq p+1$.

Choose a basis $u_1,\ldots,u_p$ of $W$ and form the matrix
\[
 M=\begin{pmatrix}v_1u_1&v_1u_2&\cdots&v_1u_p\\v_2u_1&v_2u_2&\cdots&v_2u_p\end{pmatrix}.
\]
A matrix of linear forms is \emph{$1$-generic} over a field if no invertible row and column operations over that field produce a zero entry. For $M$, this means that every nonzero linear combination of its rows has linearly independent entries. We verify this over $\overline{\kk}$. The ratio $t=v_1/v_2=(v_1w)/(v_2w)$ is a nonzero element of $\mathfrak pR_{\mathfrak p}$. If $t$ were algebraic over $\kk$, then $\kk[t]$ would be a field and $t$ would be a unit in $R_{\mathfrak p}$. Thus $t$ is transcendental. For each $\lambda\in\overline{\kk}$, let $g$ be its minimal polynomial over $\kk$. The element $g(t)$ is nonzero and homogeneous of degree zero, and hence a unit in $Q^*(R)$. In $Q^*(R)\otimes_{\kk}\overline{\kk}$ it is divisible by $t-\lambda$, so $t-\lambda$ is a unit there as well. Consequently multiplication by $av_1+bv_2=v_2(at+b)$ is injective on $W\otimes_{\kk}\overline{\kk}$ for all $[a:b]\in\mathbb P^1(\overline{\kk})$. Thus $M$ is $1$-generic over $\overline{\kk}$.

Work now over $K$ and suppress scalar extension from the notation. The map $\mathbb P^1\times\mathbb P(W)\to\mathbb P(R_1)$ sending $([a:b],[w])$ to $[(av_1+bv_2)w]$ is everywhere defined by $1$-genericity. Its image has dimension at most $p$, since $\mathbb P(W)\cong\mathbb P^{p-1}$. Each successive general hyperplane lowers the dimension of the intersection by one, and a general hyperplane misses a finite set of points. Thus $p+1$ general linear functionals on $R_1$ have no common zero on the image. Their common kernel $L$ is a general subspace of codimension $p+1$ and contains no nonzero product $(av_1+bv_2)w$. Hence every nonzero linear combination of the rows of $M$ still has linearly independent entries modulo $L$, so the reduced matrix is $1$-generic.

The polynomial ring on a basis of $R_1/L$ maps onto $R/(L)$. The entries of the reduced matrix are now linear forms in $p+1$ variables. Its $2\times2$ minors vanish in $R/(L)$ because $(v_1u_i)(v_2u_j)=(v_1u_j)(v_2u_i)$ in $R$. If $p\geq2$, these minors generate the ideal of a rational normal curve by \cite[Corollary~6.9~(b)]{E05}. Its coordinate ring $\overline{\kk}[x,y]^{(p)}$ therefore maps onto $R/(L)$. If $p=1$, the polynomial ring on $R_1/L$ itself has two variables.
\end{proof}

We next show that every nonzero graded Artinian algebra of the form $\kk[x,y]^{(p)}/J$ has a unimodal Hilbert function.

\begin{lemma}\label{lem:binary-veronese}
Let $S=\kk[x,y]$, let $p\geq1$, and let $J$ be a homogeneous ideal of $S^{(p)}$ such that $A=S^{(p)}/J$ is nonzero and Artinian. Then $S/JS$ is Artinian and $A=(S/JS)^{(p)}$. The Hilbert function $H_A(i)=\dim_{\kk}A_i$ is unimodal. For each $j\geq0$ such that $H_A(j)<pj+1$, we have $H_A(i+1)\leq H_A(i)$ for all $i\geq j$.
\end{lemma}
\begin{proof}
For each $j\geq0$ we have $[JS]_{pj}=\sum_{i=0}^j J_i S_{p(j-i)}=J_j$, since $S_{p(j-i)}$ is the degree-$(j-i)$ component of $S^{(p)}$ and the term with $i=j$ is $J_j$. Thus $A=C^{(p)}$, where $C=S/JS$. The ring $C$ is finite over $A$, since $S$ is finite over $S^{(p)}$, so $C$ is Artinian.

We may assume that $\kk$ is infinite, since extending the field preserves Hilbert functions. Set $H_C(n)=\dim_{\kk}C_n$ and let $e$ be the least nonnegative integer with $H_C(e)<e+1$. By \cite[Lemma~2.3]{I05}, $H_C(n)=n+1$ for every $0\leq n<e$, and $H_C(n+1)\leq H_C(n)$ for every $n\geq e$. The subsequence $H_A(i)=H_C(pi)$ is therefore unimodal. If $H_A(j)<pj+1$, then $pj\geq e$, so $H_A(i+1)=H_C(pi+p)\leq H_C(pi)=H_A(i)$ for every $i\geq j$.
\end{proof}

To estimate the middle entry, we compare $h_i$ and $h_{s+1-i}$.

\begin{lemma}\label{lem:opposite-entries}
Let $R$ be a $d$-dimensional standard graded Cohen--Macaulay domain over a field $\kk$ with $h$-vector $(1,h_1,\ldots,h_s)$, where $d,s\geq1$. Suppose that $[\tr_R(\omega_R)]_1\ne0$ and that $0\leq h_1-h_s\leq1$. Then
\begin{equation}\label{eq:shifted-h-comparison}
 0\leq h_i-h_{s+1-i}\leq h_1-h_s\qquad\text{for every }i\text{ with }1\leq i\leq\lceil s/2\rceil.
\end{equation}
\end{lemma}
\begin{proof}
Choose a graded fractional ideal $I\cong\omega_R(d-s)$, so that $\indeg I=0$ and $\dim_{\kk} I_0=h_s$. Lemma~\ref{lem:linear-trace-product} gives $I_0[I^{-1}]_1=[\tr_R(\omega_R)]_1\ne0$. Fix $0\ne v\in[I^{-1}]_1$. The ideal $J=vI\cong\omega_R(d-s-1)$ is a proper canonical ideal. Thus $G=R/J$ is Gorenstein of dimension $d-1$ by \cite[Proposition~3.3.18~(b)]{BH98}. Let $S$ be the polynomial ring on a basis of $G_1$ and let $Q$ be the kernel of the induced map $S\to G$. Since $J_1=vI_0$ has dimension $h_s$, we have $\dim_{\kk}G_1=d+h_1-h_s$. Hence
\[
 c:=\height Q=\dim_{\kk}G_1-\dim G=(d+h_1-h_s)-(d-1)=1+h_1-h_s.
\]
The exact sequence $0\to J\to R\to G\to0$ and graded duality give
\[
 (1-t)h_G(t)=h_R(t)-t^{s+1}h_R(t^{-1}).
\]
For $1\leq i\leq s$, the coefficient of $t^i$ on the right is $h_i-h_{s+1-i}$.

Since $c\le2$, Serre's theorem at the homogeneous maximal ideal of $S$ gives $c$ minimal generators for $Q$ \cite[Corollary~21.20]{E95}. Thus $Q$ is generated by a homogeneous regular sequence. If $h_1=h_s$, then $Q=(f)$ and $(1-t)h_G(t)=1-t^{\deg f}$. Comparing degrees gives $\deg f=s+1$, so all the indicated differences are zero.

If $h_1=h_s+1$, then $Q$ has two generators. Let $a\leq b$ be their degrees. Then
\[
 (1-t)h_G(t)=(1+t+\cdots+t^{a-1})(1-t^b),\qquad a+b=s+2.
\]
In particular, $b>\lceil s/2\rceil$. For each $i$ in the stated range, the coefficient of $t^i$ is $1$ if $i<a$ and $0$ otherwise.
\end{proof}

If $R$ satisfies $(\natur)$, then $\sqrt{[\tr_R(\omega_R)]_1R}=\m_R$. Hence $\height([\tr_R(\omega_R)]_1R)=\dim R$, and Theorem~\ref{thm:C} follows from the next result.

\begin{theorem}\label{thm:height-unimodality}
Let $R$ be a standard graded Cohen--Macaulay domain over a field $\kk$, with $h$-vector $(1,h_1,\ldots,h_s)$, where $s\geq1$ and $h_1-h_s\leq1$. Suppose that $\height([\tr_R(\omega_R)]_1R)\geq2$. Then $h_s\leq h_1$, and $h(R)$ is unimodal with $\max_i h_i-h_{\lceil s/2\rceil}\leq1$. If $\operatorname{char}\kk=0$, then the generic Artinian reduction of $R$ has the strong Lefschetz property. If $R$ is level, then $h(R)$ is concave and alternatingly increasing. In particular, in this case it is log-concave and $h_{\lceil s/2\rceil}=\max_i h_i$. If $h_1=h_s$, then $R$ is level and $h_1=\cdots=h_s$.
\end{theorem}
\begin{proof}
If $s=1$, then every linear Artinian reduction has its socle in degree one, and every nonzero linear form is a strong Lefschetz element. The other assertions follow from $h(R)=(1,h_1)$. Assume $s\geq2$ and set $d=\dim R$. Choose a graded fractional ideal $I$ isomorphic to $\omega_R(d-s)$ as a graded $R$-module. Then $\indeg I=0$ and $\dim_{\kk}I_0=h_s$ by \cite[Corollary~4.4.6~(a)]{BH98}. By Lemma~\ref{lem:linear-trace-product}, $I_0[I^{-1}]_1=[\tr_R(\omega_R)]_1$. Since $[\tr_R(\omega_R)]_1\ne0$ by the height assumption, we have $[I^{-1}]_1\ne0$. If $[I^{-1}]_1=\kk v$ for some nonzero $v$, then $[\tr_R(\omega_R)]_1R\subseteq vI$. But $vI$ is a proper canonical ideal, so it has height one by \cite[Proposition~3.3.18~(b)]{BH98}. Thus $\dim_{\kk}[I^{-1}]_1\geq2$.

Apply Lemma~\ref{lem:product-matrix} over $\kk$ with $E=[I^{-1}]_1$ and $W=I_0$. Since $\dim_{\kk}W=h_s$, we obtain a graded surjection $B=\overline{\kk}[x,y]^{(h_s)}\longrightarrow\overline R/(L)$ for every sufficiently general subspace $L\subseteq\overline R_1$ of codimension $h_s+1$, where $\overline R=R\otimes_{\kk}\overline{\kk}$. Since $h_1-h_s\leq1$, we have $\dim L=d+h_1-h_s-1\leq d$. A basis of a general such $L$ therefore extends to a linear system of parameters $\Theta$ of $\overline R$. Thus $B$ maps onto $A=\overline R/\Theta\overline R$, and $\dim_{\overline{\kk}} A_j=h_j$. Lemma~\ref{lem:binary-veronese} gives the unimodality of $h(R)$.

If $\operatorname{char}\kk=0$, Lemma~\ref{lem:binary-veronese} gives $A=C^{(h_s)}$ with $C=\overline{\kk}[x,y]/J$ Artinian. Take a strong Lefschetz element $u\in C_1$ by \cite[Proposition~4.4]{HMNW03}. For every pair $0\le a<b\le s$, multiplication by $(u^{h_s})^{b-a}$ from $A_a$ to $A_b$ is multiplication by $u^{(b-a)h_s}$ from $C_{a h_s}$ to $C_{b h_s}$. Hence $u^{h_s}$ is a strong Lefschetz element of $A$. All linear parameter reductions have Hilbert function $h(R)$, so the minors certifying maximal rank for this choice remain nonzero for generic parameters and a general linear form over their infinite coefficient field. Thus $A_{\mathrm{gen}}$ has the strong Lefschetz property.

Suppose $h_1\leq h_s$. Then $H_A(1)=h_1<h_s+1$. Lemma~\ref{lem:binary-veronese} gives $h_1\geq h_2\geq\cdots\geq h_s$, so $h_1=\cdots=h_s$. If $R$ is Gorenstein, it is level. Otherwise, $\omega_R$ has $h_s=h_1$ minimal generators in its lowest degree. The height assumption gives a nonzero $x\in[\tr_R(\omega_R)]_1$. Since $R$ is a domain, $x$ is a nonzerodivisor, so $R$ is level by \cite[Proposition~4.4~(2)]{MO}.

It remains to consider $h_1=h_s+1$. By Lemma~\ref{lem:opposite-entries}, $0\leq h_i-h_{s+1-i}\leq1$ for every $i$ with $1\leq i\leq\lceil s/2\rceil$. Hence $|h_i-h_{s+1-i}|\leq1$ for every $i$ with $1\leq i\leq s$. Since $h_1=h_s+1>h_0$, we can choose $m\geq1$ with $h_m=\max_i h_i$. The index $\lceil s/2\rceil$ lies between $m$ and $s+1-m$. Unimodality therefore gives $h_{\lceil s/2\rceil}\geq\min\{h_m,h_{s+1-m}\}\geq h_m-1$.

Now suppose that $R$ is level. Then $\overline R$ and its Artinian reduction $A$ are level as well. Every nonzero homogeneous element of $A$ has a nonzero multiple in $A_s$, since a nonzero homogeneous multiple of largest degree lies in the socle. The images of $x^{h_s}$ and $y^{h_s}$ under $B\to A$ lie in $A_1$ by the Veronese grading. For every $2\leq j\leq s$, consider the maps
\[
 A_{j-2}\xrightarrow{\ (y^{h_s},-x^{h_s})\ }A_{j-1}^{\oplus2}\xrightarrow{\ (x^{h_s},y^{h_s})\ }A_j .
\]
Their composite is zero. For every $k\geq2$, each monomial of degree $k h_s$ in $x,y$ is divisible by $x^{h_s}$ or $y^{h_s}$, since a monomial divisible by neither has degree at most $2h_s-2$. Thus $A_k=x^{h_s}A_{k-1}+y^{h_s}A_{k-1}$, so the second map is surjective. If $z\in A_{j-2}$ is in the kernel of the first map, then $x^{h_s}z=y^{h_s}z=0$. As $s-j+2\geq2$, the same equality with $k=s-j+2$ gives $z A_{s-j+2}=0$. By levelness $z=0$, so the first map is injective. Hence $h_{j-2}+h_j\leq2h_{j-1}$ for every $2\leq j\leq s$. Thus $h(R)$ is concave. For every $1\leq i<s$, concavity and the arithmetic--geometric mean inequality give $h_i\geq(h_{i-1}+h_{i+1})/2\geq\sqrt{h_{i-1}h_{i+1}}$. Squaring gives $h_i^2\geq h_{i-1}h_{i+1}$.

It remains to prove $h_i\leq h_{s-i}\leq h_{i+1}$ for every $i$ with $0\leq i<s/2$. The second inequality follows from~\eqref{eq:shifted-h-comparison} at $i+1$. If $i=0$, the first inequality is $1\leq h_s$. Suppose that $1\leq i<s/2$ and $h_i>h_{s-i}$. Then $h_k-h_{k-1}\leq-1$ for some $k$ with $i<k\leq s-i$. By concavity all subsequent differences are at most $-1$, so $h_{s-i+1}\leq h_{s-i}-1\leq h_i-2$. This contradicts~\eqref{eq:shifted-h-comparison} at $i$.
\end{proof}

\begin{example}\label{ex:middle-gap}
Consider the following two projective monomial curves, with homogeneous coordinate rings
\[
\begin{aligned}
 R&=\QQ[u^{21},u^{17}v^4,u^{12}v^9,u^9v^{12},u^8v^{13},v^{21}],\\
 S&=\QQ[u^{27},u^{22}v^5,u^{15}v^{12},u^{12}v^{15},u^{10}v^{17},v^{27}].
\end{aligned}
\]
In each ring the displayed generators have degree one. The ring $R$ is the type-four example constructed in \cite[proof of Theorem~3.11, case $d=4$]{M24}. In particular, $R$ is nearly Gorenstein and is not level. Its $h$-vector is $(1,4,7,6,3)$, with its maximum at the middle entry.

For $S$, one can check with \texttt{Macaulay2}~\cite{M2} that $S$ is Cohen--Macaulay, $\tr_S(\omega_S)=\m_S$, and $h(S)=(1,4,7,6,6,3)$. Thus $S$ is nearly Gorenstein and $h_1-h_5=1$. The middle entry specified in Theorem~\ref{thm:height-unimodality} is $h_3=6$, while $h_2=7$ is the maximum. This gives equality in the bound $\max_i h_i-h_{\lceil s/2\rceil}\le1$. Moreover, $h_3^2=36<42=h_2h_4$, so $h(S)$ is not log-concave. Thus $S$ is not level by Theorem~\ref{thm:height-unimodality}.
\end{example}

\begin{remark}\label{rem:yanagawa}
Over an algebraically closed field of characteristic zero, a standard graded Cohen--Macaulay domain with $s\geq2$ and $h_1=\cdots=h_{s-1}\geq h_s\geq2$ is level \cite[Corollary~3.11]{Y95}, \cite[Lemma~5~(2) and Proposition~6]{Y97}.
\end{remark}

\section*{Acknowledgments}
The first author was partially supported by JSPS KAKENHI Grant Number JP25K07100.
The second author was supported by JSPS KAKENHI Grant Number JP25KJ1744 and thanks Professor Sugiyama for the opportunity that led to this collaboration.
We thank Akiyoshi Tsuchiya and Ryoshun Oba for the suggestion that led to Theorem~\ref{thm:B}.

The authors used ChatGPT to examine proofs, perform computations, check references, and revise the text; they take full responsibility for all mathematical content.\par

\end{document}